\documentclass[invmat]{svjour}
\usepackage{latexsym,amssymb,amscd,amsmath}
\usepackage{graphics,color}

\def\Z{{\mathbb Z}}

\def\C{{\mathbb C}}
\def\R{{\mathbb R}}

\def\P{{\mathbb P}}
\def\CP{{\mathbb C \mathbb P}}

\newcommand{\LL}{\mathcal L}

\newcommand{\CC}{\mathcal C}
\newcommand{\II}{\mathcal I}
\newcommand{\JJ}{\mathcal J}
\newcommand{\Aa}{\mathcal A}
\newcommand{\GG}{\mathcal G}
\newcommand{\FF}{\mathcal F}
\newcommand{\MM}{\mathcal M}
\newcommand{\NN}{\mathcal N}
\newcommand{\HH}{\mathcal H}
\newcommand{\Ss}{\mathcal S}
\newcommand{\TT}{\mathcal T}
\newcommand{\OO}{\mathcal O}

\newcommand{\at}{\mathfrak A}

\newcommand{\lf}{\lfloor}
\newcommand{\rf}{\rfloor}
\newcommand{\la}{\langle}
\newcommand{\ra}{\rangle}

\newcommand{\md}{Y_{\GG}}
\newcommand{\projmd}{\overline{Y_{\GG}}} 
\newcommand{\rk}{{\rm rk}}

\newenvironment{pf}{\noindent {\bf Proof.}}{\hfill $\Box$\vspace{0.3cm}}

\begin{document}
\title{Chow rings of toric varieties defined \\ by atomic lattices}

\author{Eva Maria Feichtner\inst{1} \and Sergey Yuzvinsky\inst{2}
}                     
%
%
\institute{Department of Mathematics, ETH Zurich, 8092 Zurich, Switzerland
          \newline
          \email{feichtne@math.ethz.ch}
          \and 
          Department of Mathematics, University of Oregon, 
          Eugene OR 97403, USA
          \email{yuz@math.uoregon.edu}}
\date{Received: date / Revised version: date}
%
\maketitle
\begin{abstract}
 We study a graded algebra $D\,{=}\,D(\LL,\GG)$ over $\Z$ defined by a finite
  lattice~$\LL$ and a subset $\GG$ in $\LL$, a so-called building set.
  This algebra is a generalization of the cohomology algebras of
  hyperplane arrangement compactifications found in work of De~Concini
  and Procesi~\cite{D2}.  Our main result is a representation of~$D$,
  for an arbitrary atomic lattice~$\LL$, as the Chow ring of a smooth toric
  variety that we construct from $\LL$ and $\GG$. We describe this
  variety both by its fan and geometrically by a series of blowups and
  orbit removal. Also we find a Gr\"obner basis of the relation ideal
  of $D$ and a monomial basis of~$D$.
\end{abstract}


\section{Introduction} \label{sect_intr}
In this article we study a graded algebra $D\,{=}\,D(\LL,\GG)$ over $\Z$
that is defined by a finite lattice $\LL$ and a special subset, a so-called 
building set, $\GG$ in $\LL$.  The definition of this algebra is inspired by a
presentation for the cohomology of arrangement compactifications as it
appears in work of De~Concini and Procesi~\cite{D2}. 

In~\cite{D1,D2}
the authors studied a compactification of the complement of subspaces in
a projective space defined by a building set in the intersection
lattice $\LL$ of the subspaces. In particular they gave a description
of the cohomology algebra $H^*$ of this compactification in terms of
generators and relations.  In general, the set of defining relations
for $H^*$ is much larger than the one we propose for~$D$. However, in
the case of all subspaces being of codimension~$1$ and $\GG$ the set
of irreducibles in~$\LL$, the former can be reduced to the
latter~\cite[Prop.\ 1.1]{D2}. We show that this reduction holds for arbitrary
building sets in~$\LL$, thus giving a first geometric interpretation
of the algebra~$D(\LL,\GG)$ (compare Corollary~\ref{crl_HDP}).

Our first result about $D$ is that for an arbitrary atomic lattice
$\LL$ a larger set of relations, similar to the defining relations of
$H^*$, holds in $D$. To define the new relations for arbitrary lattices
beyond the geometric context of arrangements, we need to introduce a
special metric on the chains of $\LL$. In fact, this new set of
relations forms a Gr\"obner basis of the relation ideal which allows
us to define a basis of $D$ over $\Z$ generalizing the basis defined
in \cite{Yu} and \cite{Ga}.

Our main result about $D$ motivating its definition is Theorem
\ref{thm_ChD} which asserts that $D$ is naturally isomorphic to the
Chow ring of a smooth toric variety $X\,{=}\,X_{\Sigma(\LL,\GG)}$ 
constructed from
an atomic lattice~${\LL}$ and a building set $\GG$ in $\LL$. This
result gives a second geometric interpretation of~$D$, this time for
arbitrary atomic lattices.

We introduce the toric variety $X$ by means of its polyhedral
fan $\Sigma(\LL,\GG)$ that we build directly from $\LL$ and $\GG$.
Then we give a more geometric construction of $X$ as the result of
several toric blowups of an affine complex space and subsequent
removal of certain open torus orbits.

The article is organized as follows. In section~2, we recall the
necessary combinatorial definitions and define the algebra
$D\,{=}\,D(\LL,\GG)$. In section~3, we extend the set of relations for
$D$ to a Gr\"obner basis of the relation ideal and exhibit a basis of
the algebra. In section~4, we review the De Concini-Procesi
compactifications of arrangement complements and relate $D$ to their
cohomology algebras. Also we give some examples of the Poincar\'e
series of these compactifications using our basis. Section~5 is
devoted to the definition of the toric variety~$X$ from a pair
$(\LL,\GG)$. We prove our main theorem asserting that $D$ is naturally
isomorphic to the Chow ring of $X$. In section~6, we give another
construction of $X$ as the result of a series of toric blowups and
subsequent removal of some open orbits. Finally, in section~7, we
consider a couple of simple examples.

\bigskip


\section{The algebra $D(\LL,\GG)$} \label{sect_Dalg}

We start with defining some lattice-theoretic notions, {\em building
  sets\/} and {\em nested sets}, that provide the combinatorial
essence for our algebra definition below.  These notions, in the
special case of intersection lattices of subspace arrangements, are
crucial for the arrangement model construction of De~Concini and
Procesi~\cite{D1}.  For our purpose, we choose to present purely
order-theoretic generalizations of their notions that previously
appeared in~\cite{FK}. 

By a lattice, in this article, we mean a finite partially ordered set
all of whose subsets have a least upper bound (join, $\vee$) and a
greatest lower bound (meet, $\wedge$). The least element of any
lattice is denoted by~$\hat 0$. 
For any subset $\GG$ of a lattice $\LL$ we denote
by max$\,\GG$ the set of maximal elements of $\GG$. Also, for any
$X\,{\in}\,\LL$ we put 
$\GG_{\leq X}\,{=}\,\{G\,{\in}\,\GG\,|\,G\,{\leq}\, X\}$, 
similarly for~$\GG_{\geq}$. To denote intervals in $\LL$ we use the 
notation~$[X,Y]\,{:=}\,\{Z\,{\in}\,\LL\,|\, X\,{\leq}\,Z\,{\leq}\,Y\}$ 
for~$X,Y\,{\in}\,\LL$.

\begin{definition} \label{df_bset}
Let $\LL$ be a finite lattice. A subset $\GG$ 
in~$\LL\,{\setminus}\,\{\hat 0\}$ is called 
a {\em building set\/} in~$\LL$ if for any 
$X\,{\in}\,\LL\,{\setminus}\,\{\hat 0\}$  and
{\rm max}$\, \GG_{\leq X}{=}\{G_1,\ldots,G_k\}$ there is an isomorphism
of posets
\begin{equation*}
\varphi_X:\,\,\, \prod_{i=1}^k\,\,\, [\hat 0,G_i] \,\, 
                               \stackrel{\cong}{\longrightarrow}
                                \,\, [\hat 0,X]
\end{equation*}
with $\varphi_X(\hat 0, \ldots, G_i, \ldots, \hat 0)\, = \, G_i$ for
$i=1,\ldots, k$. We call $\max \GG_{\leq X}$ the {\em set of
factors\/} of $X$ in $\GG$.
\end{definition}

As a first easy example one can take the maximal building 
set $\LL{\setminus}\{\hat 0\}$. 
Looking at the other extreme, the elements
$X\,{\in}\,\LL\,{\setminus}\,\{\hat 0\}$ 
for which $[\hat 0,X]$ does {\em not\/}
decompose as a direct product, so-called {\em irreducibles\/} in~$\LL$, 
form the minimal building set in a given
lattice~$\LL$.

The choice of a building set $\GG$ in $\LL$ gives rise to a family of
{\em nested sets\/}. Roughly speaking these are the
subsets of $\LL$ whose antichains are
sets of factors with respect to the building set~$\GG$. The 
precise definition is as follows.

\begin{definition} \label{df_nested}
  Let $\LL$ be a finite lattice and $\GG$ a building set
  in~$\LL$.  A subset $\Ss$ in $\GG$ is called~{\em nested\/} if, for
  any set of pairwise incomparable elements $G_1,\dots,G_t$ in $\Ss$
  of cardinality at least two, the join $G_1\vee\dots\vee G_t$ 
  does not belong to $\GG$.  The nested sets in $\GG$ form an
  abstract simplicial complex, the {\em simplicial complex of nested
    sets\/} $\NN(\LL,\GG)$.
\end{definition}

For the maximal building set $\GG\,{=}\,\LL{\setminus}\{\hat 0\}$ 
the nested set complex
coincides with the order complex of $\LL{\setminus}\{\hat 0\}$.  
Smaller building sets
yield nested set complexes with fewer vertices, but allow for more
dense collections of simplexes.

An important property of a nested set is that for any two distinct 
maximal elements $X$ and $Y$ we have $X\,{\wedge}\,Y\,{=}\,\hat 0$ 
(see \cite[Prop.\ 2.5(1), 2.8(2)]{FK}).

\medskip

We now have all notions at hand to define the main character of this
article. 

\begin{definition} \label{def_D}
  Let $\LL$ be a finite lattice, $\at(\LL)$ its set of atoms, and
  $\GG$ a building set in~$\LL$. We define the 
  algebra $D(\LL,\GG)$ of $\LL$ with respect to $\GG$ as
\[ 
D(\LL,\GG) \, \, := \, \, 
             \Z\,[\{x_{G}\}_{G\in \GG}] \, \Big/ \, \II \, ,
\]
where the ideal $\II$ of relations is generated by 
\begin{eqnarray}
 \label{rel_monD}   \prod_{i=1}^t \,x_{G_i} & & \qquad 
          \mbox{for }\,\,\{G_1,\ldots ,G_t\}\,\not\in \,\NN(\LL,\GG) \,,
\end{eqnarray}
and 
\begin{eqnarray}
 \label{rel_linD}   \sum_{G\geq H}\, x_G & & \qquad  
                            \mbox{for }\, H \in \at(\LL)\, .  
\end{eqnarray}
\end{definition} 

Note that the algebra $D(\LL,\GG)$ is a quotient of the face
algebra of the simplicial complex $\NN(\LL,\GG)$.
Although $D$ is defined for an arbitrary lattice our main constructions and
results make sense only for atomic lattices, i.e., lattices in which any 
element is the join of some atoms. Thus we will restrict our
considerations to this case.

In the special case of $\LL$ being the intersection lattice of an
arrangement of complex linear hyperplanes and $\GG$ being the minimal
building set in~$\LL$, this algebra appears in work of De~Concini and 
Procesi~\cite{D2}. It is the cohomology algebra of a compactification
of the projectivized arrangement complement; for details we refer to
section~\ref{sect_arrgtcomp}.


\section{Gr\"obner basis} \label{sect_grobner}
 
The set of generators of the ideal $\II$ in Definition~\ref{def_D}, while 
being elegant, is too small for being a Gr\"obner basis of this ideal. In this
section, we extend this set to a Gr\"obner
basis. In particular, we will obtain a $\Z$-basis of $D(\LL,\GG)$.

To define the larger set of relations 
we need to introduce a metric on chains in~$\LL$.
 
\begin{definition}\label{def_d}
Let $\LL$ be an atomic lattice and $X,Y\,{\in}\,\LL$ with $X\,{\leq}\,Y$.
We denote by $d(X,Y)$ the minimal number of atoms $H_1,\ldots,H_d$ in~$\LL$
such that $Y\,{=}\,X\,{\vee}\,\bigvee_{i=1}^d H_i$.
\end{definition}

\noindent
The following four properties of the function $d$ are immediate:
\begin{quote}
\begin{itemize}
\item[(i)] 
$d(X,Z)\,{\geq}\,d(Y,Z)$ for $X,Y,Z\,{\in}\,\LL$ with 
$X\,{\leq}\,Y\,{\leq}\,Z$.
Notice that equality is possible even if all three $X,Y,$ and $Z$ are
distinct. Also it is not necessarily true that $d(X,Y)\leq d(X,Z)$.

\item[(ii)] $d(X,Y){+}d(Y,Z)\,{\geq}\,d(X,Z)$ for $X,Y,Z\,{\in}\,\LL$ with 
$X\,{\leq}\,Y\,{\leq}\,Z$. 

\item[(iii)] $d(X\vee Z,Y\vee Z)\,{\leq}\,d(X,Y)$ for $X\,{\leq}\,Y\,{\in}\,\LL$
and $Z\,{\in}\,\LL$ arbitrary.

\item[(iv)] 
$d(A,A\vee B)\,{\leq}\,d(A\wedge B,B)$ for $A,B\,{\in}\,\LL$\,. 
\end{itemize}
\end{quote}

\noindent
For example, (iv) follows from the fact that if $(A\wedge B)\vee
\bigvee_i H_i\,{=}\,B$ for some atoms $H_1,\ldots, H_d$ then
$A\vee\bigvee_iH_i\,{=}\,A\vee B$. If $\LL$ is geometric (for instance, the
intersection lattice of a hyperplane arrangement) then $d(X,Y)\,{=}\,\rk
Y\,{-}\,\rk X$ whence in (ii) equality holds and (iv) is the semimodular 
inequality.

Now we can introduce the new set of generators for $\II$. The new
relations are analogous to the defining relations for the cohomology
algebra of the compactification of the complement of an arrangement of
projective subspaces described in~\cite{D1}.

\begin{theorem} \label{prop_Hmd} 
The ideal of relations~$\II$ in Definition \ref{def_D}
is generated by polynomials of the following type:
\begin{eqnarray}
 \label{rel_monH} && h_{\Ss}= \prod_{G\in \Ss} \,x_{G} \hspace{1.5cm}
                     \mbox{for }\,\, \Ss\not\in \NN(\LL,\GG)\, , \\
 \label{rel_monlinH}&& g_{\HH,B}=\prod_{i=1}^k \,x_{A_i} \,  
                          \Big(\sum_{G\geq B}\, x_G \,\Big)^d\, , 
\end{eqnarray}
where $A_1,\ldots, A_k$ are maximal elements in a nested set $\HH\,{\in}\,
\NN(\LL,\GG)$, $B\,{\in}\,\GG$ with $B\,{>}\,A=\bigvee_{i=1}^k A_i$, and
$d\,{=}\,d(A,B)$. 
\end{theorem}

\begin{pf}
First notice that polynomials (\ref{rel_monD}) and (\ref{rel_linD}) are among 
polynomials $h_{\Ss}$ and $g_{\HH,B}$.
   (To see that
  polynomials (\ref{rel_linD})  are among $g_{\HH,B}$ choose
  $\HH=\emptyset$, and $B\,{=}\,H\,{\in}\,\at(\LL)$. Here and everywhere we use
the usual agreement that the join of the empty set is $\hat 0$.)
  Hence it is left to show that any $g_{\HH,B}$ is in $\II$, i.e., it is a
combination of polynomials
~(\ref{rel_monD}) and~(\ref{rel_linD}).

We prove our claim by induction on $d$. \newline
{\underline{$d=1$.}} 
Choose an atom $H$ of $\LL$ with $H \vee A\,{=}\, B$. Then using
(\ref{rel_linD}) we have
\begin{eqnarray} \label{eqn_d=1}
    \prod_{i=1}^k\,\, x_{A_i} \, \, 
\Big(\,\sum_{G\geq H}\, x_G \,\Big) \, \in\II.
\end{eqnarray}
We want to show that for any $G\geq H$, $\{G,A_1,\ldots,
A_k\}\,{\in}\, \NN\,{=}\, \NN(\LL,\GG)$ implies that $G\geq B$. Then, any summand
with $G\not\geq B$ can be omitted from~(\ref{eqn_d=1}) using
polynomials~(\ref{rel_monD}), and we obtain $g_{\HH,B}\in\II$ for $d\,{=}\,1$.

First note that $G$ cannot be smaller than or equal to any of the $A_i$,
$i=1,\ldots,k$, since $G\leq A_i$ would imply $H\leq A_i$ contradicting
the choice of~$H$.

Assume that $G$ is incomparable with $A_1,\ldots, A_s$ for some
$s\,{\geq}\, 1$, and $G\,{\geq}\,A_i$ for $i\,{=}\,s{+}1,\ldots, k$.
Since $\{G,A_1,\ldots, A_k\}\,{\in}\, \NN$ these elements 
are the factors of the \mbox{$\GG$-de}composition in
\[
\widetilde G\,\, := \,\, G \vee \bigvee_{i=1}^s A_i  
        \,\, = \,\, G \vee \bigvee_{i=1}^k A_i  
        \,\, \geq \,\, H \vee \bigvee_{i=1}^k A_i \, \, = \, \, B\, .
\] 
Since $B\,{\in}\,\GG$, the elements $A_i$, $i=1,\ldots ,s$, are not 
maximal in $\GG$
below $\widetilde G$, which contradicts the $A_i$ being factors of
$\widetilde G$.

We conclude that $G$ is comparable with, i.e., larger than all $A_i$
whence $G\geq \bigvee_{i=1}^k A_i \, \vee \, H \, = \, B$.

\smallskip
\noindent
{\underline{$d>1$.}} Choose an atom $H$ of $\LL$ from the set of atoms
in the definition of $d(A,B)$. Then
$A\,{<}\,A\vee H\,{<}\,B$.  Using~(\ref{rel_linD}) we have
\begin{eqnarray} \label{eqn_d>1}
    \prod_{i=1}^k\,\, x_{A_i} \, \, 
              \Big(\,\sum_{G\geq H}\, x_G \,\Big) \,\,
           \Big(\,\sum_{G\geq B}\, x_G \,\Big)^{d-1}                          
                                                        \, \in\II .
\end{eqnarray}
We show, using polynomials ~(\ref{rel_monD})
and~(\ref{rel_linD}) and the induction hypothesis, that any $G$ with
$G\,{\not\geq}\,B$ can be omitted from the first sum modulo~$\II$.

Let $G_0\,{\in}\,\GG$, $G_0\,{\geq}\,H$ but $G_0\,{\not\geq}\,B$. 
Using polynomials
~(\ref{rel_monD}) we can assume that
$\{G_0,A_1,\ldots,A_k\}\,{\in}\,\NN$. Due to the choice of $H$,
$G_0$ cannot be smaller than any of the $A_i$. Further note that if
$G_0$ is incomparable with say $A_1,\ldots A_s$, $s\leq k$, then it is
incomparable also with
all $A_1,\ldots, A_k$. Indeed the join $G_0\vee A_1\vee \ldots \vee
A_s\,{=}\,G_0\vee A_1\vee \ldots \vee A_k$ is a $\GG$-decomposition.
Hence the two following cases remain to be considered.

\noindent
{\bf Case 1.} {\em $G_0$ is comparable with all $A_i$, $i=1,\ldots,k$,
  hence $G_0\geq A$.} 
\newline Our goal is to rewrite
\begin{eqnarray} \label{eqn_case1}
           x_{G_0}\, \, \Big(\,\sum_{G\geq B}\, x_G \,\Big)^{d-1}
\end{eqnarray}
modulo $\II$
so that it contains an expression of the form~(\ref{rel_monlinH}) with
exponent $<d$ as a factor. First observe that
$G_0{\vee}B\,{\in}\,\GG$ since $G_0,B\,{\in}\,\GG$ but $H< G_0 \wedge
B$~\cite[Thm.~2.3, 3b']{D1}.  The building set element $G_0{\vee}B$ is to
take the role of $B$ in~(\ref{rel_monlinH}).

Let $G\in \GG$ with $G\geq B$. We want to show that any $G$ with $G
\not\geq G_0\vee B$ can be omitted from~(\ref{eqn_case1}) modulo $\II$.  We can
assume that $\{G,G_0\}\in \NN$. If $G\leq G_0$ then $B\leq G_0$,
contradicting the choice of~$G_0$. If $G$ and $G_0$ were incomparable then
$G\vee G_0\not \in \GG$ contradicting the fact that they both are greater than
$H$.  Hence $G\geq G_0$ and thus $G\geq G_0{\vee}B$.

Thus (\ref{eqn_case1}) reduces to
\begin{eqnarray}\label{eqn_case1D} 
           x_{G_0}\, \, \Big(\,\sum_{G\geq G_0\vee B}\, x_G \,\Big)^{d-1}\, .
\end{eqnarray}
Using properties (iv) and (i) of our metric $d$ we obtain 
\begin{eqnarray}
d( G_0,G_0\vee B)  \leq  d(G_0 \wedge B,B)  \leq  
d(A \vee H,B) < d\, .
\end{eqnarray}
Hence (\ref{eqn_case1D}) contains a polynomial of the
form~(\ref{rel_monlinH}) with exponent $<d$ as a factor whence it
lies in $\II$ by induction hypothesis.

\noindent
{\bf Case 2.} {\em $G_0$ is incomparable with $A_1,\ldots,A_k$.}
\newline Since $\{G_0,A_1,\ldots A_k\}\in \NN$ we have $\widetilde G_0:=
G_0\vee A_1 \vee \ldots \vee A_k \not \in \GG$.  We want to rewrite
\begin{eqnarray} \label{eqn_case2}
     \Big(\,\prod_{i=1}^k\,\, x_{A_i}\,\Big) \,\, x_{G_0}\, \, 
                            \Big(\,\sum_{G\geq B}\, x_G \,\Big)^{d-1}
\end{eqnarray}
modulo $\II$
so that it contains a polynomial of the form~(\ref{rel_monlinH}) with
exponent $<d$ as a factor.

Observe that $\widetilde G_0\vee B = G_0\vee B$, and, as in Case~1,
$G_0\vee B \in \GG$. This time, $\widetilde G_0\vee B= G_0\vee B$ is
to take the role of~$B$, and $\widetilde G_0$ the role of $A$
in~(\ref{rel_monlinH}).

As in Case~1, we see that
\begin{eqnarray*} 
\lefteqn{
     \Big(\, \prod_{i=1}^k\,\, x_{A_i}\,\Big) \,\, x_{G_0}\, \, 
                            \Big(\,\sum_{G\geq B}\, x_G \,\Big)^{d-1} 
     \quad \equiv } \\
& & \qquad \qquad \quad
      \Big(\, \prod_{i=1}^k\,\, x_{A_i}\, \Big)\,\, x_{G_0}\, \, 
           \Big(\,\sum_{G\geq G_0{\vee} B}\, x_G \,\Big)^{d-1}\, {\rm modulo}\
\II,
\end{eqnarray*}
arguing as before for nested pairs $\{G,G_0\}$. 

Now the right hand side has a factor of the 
form~(\ref{rel_monlinH}) with exponent $<d$ because again
\begin{eqnarray*}
d(\widetilde G_0,\widetilde G_0\vee B) \leq d(\widetilde G_0 \wedge B,B)
 \leq  d(B,A \vee H)  <  d\, .
\end{eqnarray*}
 This implies that the right hand side lies in $\II$ by induction hypothesis
which completes the proof.
\end{pf} 

The main feature of the new generating set is that it is a Gr\"obner basis of
$\II$. As the main reference for Gr\"obner bases we use
\cite{Ei}. Fix a linear order on $\GG$ that refines the reverse of
the partial order on $\LL$. It defines a  lexicographic order on the monomials
which we use in the following theorem.

\begin{theorem} \label{thm_Gr}
The generating system (\ref{rel_monH}) and (\ref{rel_monlinH}) is a
Gr\"obner basis of~$\II$.
\end{theorem}

\begin{pf}
  To prove that a set of monic polynomials is a Gr\"obner basis for
  the ideal it generates it suffices to consider all pairs of their
  initial monomials with a common indeterminate, compute their
  syzygies, and show that these syzygies have standard expressions in
  generators (without remainders). We will prove this by a
  straightforward calculation. To make the calculation easier to
  follow we will use several agreements. For any polynomial $p\in\II$
  we will be dealing with, we will exhibit a generator $g$ whose
  initial monomial $in(g)$ divides a monomial $\mu$ of $p$ and call
  $p-c(\mu){\frac{\mu}{in(g)}}g$ the reduction of $p$ by $g$ (here
  $c(\mu)$ is the coefficient of $\mu$ in $p$).  Reducing a polynomial
  all the way to~$0$ gives a standard expression for it.  Also since
  reduction by monomial generators is very simple we will not name
  specific generators of the form $h_{\Ss}$ but just call this
  reduction $h$-equivalence.

We use certain new notation in the proof. For each $\Ss\,{\subset}\,\GG$ put
$\pi_{\Ss}\,{=}\,\prod_{A\in \Ss}x_A$ and for any $B\,{\in}\,\GG$ put 
$y_B=\sum_{Y\in\GG_{>B}}x_Y$.

Now we consider pairs $(g_1,g_2)$ of generators of $\II$ of several types.

{\bf 1.} 
At least one of the generators is $h_{\Ss}$. If they both are of
this type then the syzygy is~$0$. If the other one is $g_{\HH,B}$ with
$B\not\in {\Ss}$ then the syzygy is divisible by $h_{\Ss}$ whence
$h$-equivalent to 0. Finally if $B{\in}{\Ss}$ then the only nontrivial 
case is where 
$T\,{=}\,(\Ss\cup\HH){\setminus}\{B\}\,{\in}\,\NN{=}\NN(\GG,\LL)$. 
Notice that
then $\Ss\cup\HH\not\in\NN$. The syzygy is $h$-equivalent to
$\pi_Ty_B^{d(A,B)}$ where $A\,{=}\,\bigvee_{X\in\HH}X$ as usual. Put $\bar
A\,{=}\,\bigvee_{X\in T}X$.  If $X\,{\in}\,\GG_{>B}$ and 
$X\,{\leq}\, \bar A\vee B$ then
$X$ cannot form a nested set with $T$. Indeed, if it did then $\bar
A\vee B =X\vee\bar A\not\in\GG$ contradicting
$T\cup\{X\}\not\in\NN$. Similarly, if $X\in\GG_{>B}$ and $X$ is
incomparable with $\bar A$ then $X$ cannot form a nested set with $T$.
Indeed if they did then $X\vee(\bar A\vee B)=X\vee \bar A\not\in \NN$
implying that $X$ forms a nested set with $\Ss\,{\cup}\,\HH$. This would
contradict $X\,{>}\,B$.

Now using property (i) of the metric $d$ we can reduce the syzygy to~$0$ 
by $g_{T,\bar A\vee B}$.

For the rest of the proof we need to consider only pairs with
$g_i=g_{\HH_i,B_i}$ ($i=1,2)$. We denote the exponent of $x_{B_i}+y_{B_i}$ 
in $g_i$ by $d_i$.

{\bf 2.}  Suppose $B_1\,{\not=}\,B_2$ and $B_i\,{\not\in}\,\HH_j$.
In this case the syzygy is
$$\pi_{\HH_2\setminus \HH_1}g_1(g_1-in(g_1))-\pi_{\HH_1\setminus\HH_2}g_2(g_2
-in(g_2))$$
and this is in fact a standard expression for it. (Here and to the end of the
proof we use $\pi_{\Ss}$ for arbitrary subsets $\Ss$ of $\LL$ meaning that 
if $\Ss$ is not nested the product is $h$-equivalent to~$0$.)

{\bf 3.} Suppose $B_1\,{=}\,B_2\,{=}\,B$ and $d\,{=}\,d_2-d_2\,{\geq}\, 0$.
Then the syzygy is 
$$\pi_{\HH_1\cup\HH_2}\,[\,x_B^d(x_B+y_B)^{d_1}-(x_B+y_B)^{d_2}\,]$$
and it reduces to~$0$ by $g_1$.

{\bf 4.}  At last, suppose $B_1 \in \HH_2$. 
Put $\HH = (\HH_1\,{\cup}\,\HH_2)\,{\setminus}\,\{B_1\}$
and $x_{B_i}\,{=}\,x_i$, $y_{B_i}\,{=}\,y_i$.
Then the syzygy is
$$s=\pi_{\HH}\,[\,(x_1+y_1)^{d_1}x_2^{d_2}-x_1^{d_1}(x_2+y_2)^{d_2}\,]\,.$$
Adding to $s$ the polynomial $f=\pi_{\HH}(x_1+y_1)^{d_1}[(x_2+y_2)^{d_2}-x_2^{d_2}]$
we obtain
$$s'=s+f=\pi_{\HH}\,[\,(x_1+y_1)^{d_1}-x_1^{d_1}\,]\,(x_2+y_2)^{d_2}\,.$$
Notice that $f$ is divisible by $g_1$ and $in(f)\leq in(s)$.
Thus it suffices to reduce $s'$ to~$0$.
Also by $g_2$ we can immediately reduce $s'$ to 
$$s''=\pi_{\HH}\,y_1^{d_1}(x_2+y_2)^{d_2}\,.$$

For the next steps we sort out summands of $y_1$.
Using property (i) of the metric $d$ we can delete the summands $x_Y$ with
$B_1<Y<B_2$ reducing by $g_{\HH\cup\{Y\},B_2}$. The sum of all summands
$x_Y$ with $Y\geq B_2$ forms $\pi_{\HH}(x_2+y_2)^{d_1+d_2}$
that reduces to 0 by $g_{\HH,B_2}$. Indeed, denote the join 
of $\HH_i$ by $C_i$ and the join of $\HH_2\,{\setminus}\,\{B_1\}$ 
by $C_2'$. This
gives the join of $\HH$ as $C_1\vee C_2'$. Then,
using properties (ii) and (iii) of the metric~$d$, we have
\begin{eqnarray*}
d(C_1\vee C_2',B_2)& \leq & d(C_1\vee C_2',B_1\vee C_2')\, +\, 
                            d(B_1\vee C_2',B_2)\\
                   & \leq & d(C_1,B_1)\, +\, d(C_2,B_2)\, \, = \,\, 
                            d_1+ d_2\, .
\end{eqnarray*}

After the reductions in the previous paragraph we are left with
a sum each summand of which is divisible by a polynomial
$$t_Z=\pi_{\HH}x_Z(x_2+y_2)^{d_2}\, , $$
where $Z\,{\in}\,\GG_{>B_1}$, $Z$ is incomparable with $B_2$, and 
$\HH\,{\cup}\,\{Z\}\,{\in}\,\NN$.
To reduce this polynomial we sort out the summands in the second sum.
If $Y\in\GG_{\geq B_2}$ is not greater than or equal to $Z\vee B_2$ then it is
incomparable with $Z$ whence $\{Z,Y\}\,{\not\in}\,\NN$ 
since $B_1\,{<}\,Z,Y$. This implies
that $t_Z$ is $h$-equivalent to 
$$t'_Z\,=\,\pi_{\HH}x_Z(\sum_{Y\geq B_2\vee Z}x_Y)^{d_2}. $$
Finally $t'_Z$ reduces to~$0$ by $g_{\HH\cup\{Z\},B_2\vee Z}$ 
since, by property (iii) of the metric $d$, we have
$$d(C_2'\vee Z,B_2\vee Z)\,=\,d(C_2\vee Z,B_2\vee Z)\,\leq\, d(C_2,B_2)
\,=\,d_2.$$
This reduction completes the proof.
\end{pf} 

\begin{corollary} \label{crl_basis}
The following monomials form a $\Z$-basis of the algebra $D(\LL,\GG)$:
$$\prod_{A\in\Ss} x_A^{m(A)}\, ,$$
where $\Ss$ is running over all nested subsets of $\GG$ and 
$m(A)\,{<}\,d(A',A)$, $A'$ being the join of $\Ss\,{\cap}\,\LL_{<A}$.
\end{corollary}

If $\LL$ is the intersection lattice of a complex central hyperplane
arrangement then this basis coincides with the basis exhibited in \cite {Yu}.
In the next section we will give some examples of computing the Hilbert series
of the algebra using this basis.


\section{Arrangement Compactifications} \label{sect_arrgtcomp}

As we mentioned before, for a geometric lattice the metric $d$ defined
in section~\ref{sect_grobner} coincides with the difference of ranks. 
This holds in
particular for the intersection lattice of a hyperplane arrangement.
In this setting and for $\GG$ being the minimal building set, the
algebra $D(\LL,\GG)$ appeared in \cite{D2} as the cohomology algebra
of a compactification of the projectivized arrangement complement.
From our work in previous sections we can conclude that for {\em any\/}
building set $\GG$ in $\LL$ the algebra $D(\LL,\GG)$ can be
interpreted geometrically as the cohomology algebra of the
corresponding arrangement compactification.

We first review the construction of arrangement models due 
to De~Con\-cini and 
Procesi in the special case of complex hyperplane arrangements~\cite{D1}. 

Let $\Aa\,{=}\,\{H_1,\ldots,H_n\}$ be an arrangement of complex linear
hyperplanes in~$\C^d$. Factoring by $\bigcap H_i$ if needed, we can
assume $\Aa$ to be essential, i.e., $\bigcap
H_i\,{=}\,\{0\}$.
 The combinatorial data of such an
arrangement is customarily recorded by its intersection
lattice~$\LL(\Aa)$, i.e., the poset of intersections of all subsets
of hyperplanes ordered
by reverse inclusion. The greatest element of $\LL(\Aa)$ is 0 and the least
element is  $\C^d$. Let $\GG\,{\subseteq}\,\LL(\Aa)$ be a building
set in $\LL(\Aa)$, and let us assume here that $0\,{\in}\,\GG$.

We define a map on $\MM(\Aa)\,{:=}\,
\C^d\,{\setminus}\,\bigcup \Aa$, the complement of the arrangement,
\[
    \Phi: \quad \MM(\Aa) \, \, \longrightarrow \, \, 
                    \C^d\, \times \, \prod_{G\in\GG}\, \P(\C^d/G)\, ,
\]
where $\Phi$ is the natural inclusion into the first factor and
the natural projection to the other factors restricted to $\MM(\Aa)$.
The map $\Phi$ defines an embedding of $\MM(\Aa)$ in the right hand side space
and we let $\md$ denote
the closure of its image. The space $\md$ is a smooth algebraic variety
containing $\MM(\Aa$) as an open set. The complement
$\md\,{\setminus}\,\MM(\Aa)$ is a divisor with normal crossings with
irreducible components indexed by building set elements. An
intersection of several components is non-empty (moreover, transversal and
irreducible) if and only if the index set is nested as a subset 
of~$\GG$~\cite[3.1,3.2]{D1}.

There is a projective analogue of $\md$.
Consider the projectivization~$\P\Aa$ of~$\Aa$, i.e., the family of
codim~$1$ projective spaces $\P H$ in $\CP^{d-1}$ for $H\,{\in}\,\Aa$. The
following construction yields a compactification of the complement
$\MM(\P\Aa)\,{:=}\,\CP^{d-1}\,{\setminus}\, \bigcup \P \Aa$. The map
$\Phi$ described above is $\C^*$-equivariant, where $\C^*$ acts by
scalar multiplication on $\MM(\Aa)$ and on $\C^d$, and trivially on
$\prod_{G\in\GG}\, \P(\C^d/G)$. We obtain a map
\[
  \overline \Phi: \quad \MM(\P \Aa) \, \, \longrightarrow \, \, 
                    \CP^{d-1}\, \times \, \prod_{G\in\GG}\, \P(\C^d/G)\, ,
\]
and again take the closure of its image to define a model $\projmd$
for $\MM(\P \Aa)$. The space $\projmd$ is a smooth projective variety and
the complement $\projmd\,{\setminus}\,\MM(\P \Aa)$ is a divisor
with normal crossings. Irreducible components are indexed by building
set elements in $\GG^0:=\GG\,{\setminus}\,\{\{0\}\}$, and
intersections of irreducible components are non-empty if and only if
corresponding index sets are nested in~$\GG$.

Geometrically, the arrangement models $\md$ and $\projmd$ are related
as follows. The model $\md$ is the total space of a line bundle over $\projmd$;
in fact, it is the pullback of the tautological bundle on $\CP^{d-1}$
along the canonical map $\projmd\,{\rightarrow}\,\CP^{d-1}$. In
particular, $\projmd$ is isomorphic to the divisor in $\md$ associated 
to $0$~\cite[4.1]{D1}.

\begin{example} \label{ex_braid}
Let $\Aa_{n-1}$ denote the rank $n{-}1$ complex {\em braid
arrangement}, i.e., the family of partial diagonals, $H_{i,j}:
z_j{-}z_i\,{=}\,0$, $1{\leq}i{<}j{\leq}n$, in $\C^n$. Its intersection
lattice~$\LL(\Aa_{n-1})$ equals the lattice $\Pi_n$ consisting of the
set partitions of~$[n]:=\{1,\ldots,n\}$ ordered by reverse
refinement.  The set $\FF$ of partitions with
exactly one block of size $\geq 2$ forms 
the minimal building set in~$\Pi_n$.  
The De~Concini-Procesi
arrangement compactification $\overline{Y_{\FF}}$ is isomorphic to the
Deligne-Knudson-Mumford compactification of the moduli space
$M_{0,n+1}$ of $n{+}1$-punctured complex projective
lines~\cite[4.3]{D1}.  
\end{example}

In the more general setting of affine models for complex subspace
arrangements, De~Concini and Procesi provide explicit presentations
for the cohomology algebras of irreducible components of divisors and
of their intersections in terms of generators and
relations~\cite[\S5]{D1}. As mentioned above, the compactification of
a complex hyperplane arrangement $\projmd$ is isomorphic to the
divisor associated with the maximal building set element in the
corresponding affine model. We recall a description of its integral
cohomology algebra.

\begin{proposition} \label{prop_Hmd1} 
  {\rm (\cite[Thm.\ 5.2]{D1})}
  Let $\Aa$ be an essential arrangement of complex hyperplanes,
  $\LL\,{=}\,\LL(\Aa)$ its intersection lattice, and $\GG$ a building set in
  $\LL$ containing~$\{0\}$. Then the integral cohomology algebra
  of the arrangement compactification $\projmd$ can be described as 
\[ 
 H^*(\projmd) \, \, \cong \, \, 
          \Z\,[\{c_{G}\}_{G\in \GG}] \, \Big/ \, \JJ \, ,
\]
with generators $c_G$, $G\,{\in}\,\GG$, corresponding to the
cohomology classes of irreducible components of the normal crossing
divisor, thus having degree~$2$. \newline
The ideal of relations $\JJ$ is
generated by polynomials of the following type:
\begin{eqnarray}
  & & \prod_{i=1}^t \,c_{G_i} \hspace{1cm}
                     \mbox{for }\,\,\{G_1,\ldots ,G_t\}
                       \not\in \NN(\LL,\GG)\, ,  \\
  & & \prod_{i=1}^k \,c_{A_i} \,  
                          \Big(\sum_{G\geq B}\, c_G \,\Big)^d\, , 
\end{eqnarray}
where $A_1,\ldots, A_k$ are maximal elements in a nested set $\HH\,{\in}\,
\NN(\LL,\GG)$, $B\,{\in}\,\GG$ with $B\,{>}\,\bigvee_{i=1}^k A_i$, and
$d\,{=}\, {\rm codim}_{\C} B -{\rm codim}_{\C}\!\bigvee_{i=1}^k A_i$.
\end{proposition}

Comparing Proposition \ref{prop_Hmd1} with Theorem \ref{prop_Hmd}, we have 
a generalization of Proposition~1.1 from \cite{D2}, where only the case of 
$\GG$ being the minimal building set, i.e., the set of irreducibles, is 
considered.

\begin{corollary} \label{crl_HDP}
  Let $\Aa$ be an essential arrangement of complex hyperplanes,
  $\LL\,{=}\,\LL(\Aa)$ its intersection lattice, and $\GG$ a 
  building set in~$\LL$ containing~$\{0\}$. Then the cohomology 
  algebra of the 
  arrangement compactification $\projmd$ is isomorphic to the 
  algebra~$D(\LL,\GG)$ defined in section~\ref{sect_Dalg}:
\[
        H^*(\projmd)\, \, \cong \,\, D(\LL,\GG)\, .
\]
\end{corollary}

\medskip

In the rest of the section we will give several examples of the Poincar\'e 
series
for compactifications of hyperplane arrangement
complements. This means we compute the Hilbert series of $D(\LL,\GG)$.
We restrict our computations to the compactifications with
$\GG$ being the maximal building set $\LL\setminus\{\hat 0\}$, 
although they can be easily generalized to arbitrary $\GG$.

For these examples we use the basis of 
$D(\LL)\,{=}\,D(\LL,\LL\,{\setminus}\,\{\hat 0\})$ 
from Corollary~\ref{crl_basis}.
In the considered case the basic monomials are paramet\-rized 
by certain flags in
$\LL\,{\setminus}\,\{\hat 0\}$ with multiplicity assigned to their 
elements. The upper bounds
for multiplicities allow us to write the Hilbert series of ${\rm
D}(\LL)$ in the following form. For each sequence
$r$ of natural numbers, $r\,{=}\,(0{=}r_0\,{<}\,r_1\,{<}\,\cdots\,$ 
${<}\,r_k\,{\leq}\,\rk\LL)$ 
denote by $f_{\LL}(r)$ the number
of flags in $\LL$ whose sequence of ranks equals $r$. Set $k\,{=}\,k(r)$
and call it the length of $r$.  Then we have 
$$H({\rm
D}(\LL),t)=1+\sum_r\Biggl[\prod_{i=1}^{k(r)}{\frac{t(1-t)^{r_i-r_{i-1}-1}}
{1-t}}\Biggr]f_{\LL}(r).$$ 
Here, $r$ runs over all sequences as above and we
use the agreement ${\frac{t(t-1)^0}{t-1}}=1$.

In some important cases one can give more explicit descriptions of the
numbers $f_{\LL}(r)$ whence of the Hilbert series. We consider two
such cases.

\smallskip
\noindent
{\bf Generic arrangements.} For arrangements from this class, the
intersection lattice~$\LL$ is defined by the number $n$ of atoms
and the rank~$\ell$. We use both pieces of notation: $\LL$ and
$\LL(n,\ell)$.  The number of elements of $\LL$ of rank $\ell'\,{<}\,\ell$
is $\binom{n}{\ell'}$ and for every $X\,{\in}\,\LL$ of rank~$\ell'$ the
lattice $\{Y\,{\in}\,\LL\,|\,Y\,{\geq}\, X\}$ is isomorphic to
$\LL(n{-}\ell',\ell{-}\ell')$. This immediately implies the following
formula:
$$f_{\LL}(r)=\prod_{i=1}^{k}\binom{n-r_{i-1}}{r_i-r_{i-1}}\, ,$$ where
$k\,{=}\,k(r)$ if $r_{k(r)}\,{<}\,\ell$ and $k\,{=}\,k(r){-}1$ otherwise.  
This gives
\begin{eqnarray*}
\lefteqn{
H({\rm D}(\LL(n,\ell)),t)\, \, = } \\
& & 1 \,+\, \sum_r\Biggl\{\Biggl[1+{\frac{t(1-t)^{\ell-r_k-1}}{1-t}}
\Biggr]
\prod_{i=1}^{k(r)}{\frac{t(1-t)^{r_i-r_{i-1}-1}}{1-t}}
{\binom{n-r_{i-1}}{r_i-r_{i-1}}}\Biggr\}\, ,
\end{eqnarray*}
where the summation now is over all $r$ with the extra condition
$r_{k(r)}\,{<}\,\ell$ and we again use the agreement 
${\frac{t(t-1)^0}{t-1}}=1$.

\smallskip
\noindent
{\bf Braid arrangements.} For the rank $n{-}1$ complex braid arrangement 
(compare
Example~\ref{ex_braid}) the intersection lattice is given by the partition 
lattice 
$\Pi_n$ of set partitions of $[n]\,{:=}\,\{1,\ldots,n\}$ ordered by reverse 
refinement. Observe that the rank of a partition $\pi$ coincides with 
$n{-}|\pi|$
where $|\pi|$ is the number of blocks of the partition. Thus the
number of elements of $\Pi_n$ of rank $\ell$ is $p_{n-\ell}(n)$ that
is the number of partitions of $[n]$ in $n{-}\ell$ blocks.  For every
$X\,{\in}\,\Pi_n$ of rank $\ell$ the lattice 
$\{Y\,{\in}\,\Pi_n\,|\,Y\,{\geq}\, X\}$ is
isomorphic to $\Pi_{n-\ell}$. This immediately implies the following
formulas: 
$$f_{\Pi_n}(r)=\prod_{i=1}^{k(r)}p_{n-r_{i}}(n-r_{i-1})$$
and 
$$H({\rm D}(\Pi_n),t)=1+\sum_r\Biggl[\prod_{i=1}^{k(r)}
{\frac{t(1-t)^{r_i-r_{i-1}-1}}{1-t}}p_{n-r_{i}}(n-r_{i-1})\Biggr]\, ,$$
where the summation is over all $r$.


\section{The toric variety $X_{\Sigma(\LL,\GG)}$} \label{sect_torv}

In this section we present another geometric interpretation
of the algebra~$D(\LL,\GG)$, this time for an arbitrary atomic lattice
$\LL$. For a given building set~$\GG$ 
in~$\LL$ we construct a toric variety~$X_{\Sigma(\LL,\GG)}$ and 
show that its Chow ring is isomorphic to the 
algebra~$D(\LL,\GG)$.

Given a finite lattice $\LL$ with set of atoms
$\at(\LL)\,{=}\,\{A_1,\ldots, A_n\}$, we will frequently use the
following notation: For $X\,{\in}\,\LL$, denote the set of
atoms below~$X$ by $\lf X
\rf\,{:=}\,\{A\,{\in}\,\at(\LL)\, |\, X\,{\geq}\,A\}$. 
 Define characteristic
vectors $v_X$ in $\R^n$ for $X\,{\in}\,\LL$ with coordinates
\[
(v_X)_i \, \, := \, \, \left\{ 
\begin{array}{ll}
1 & \mbox{ if }\, A_i \in \lf X \rf, \\
0 & \mbox{ otherwise}, \qquad \,\,\, \mbox{ for }\, i=1,\ldots,n.
\end{array}
\right. 
\]
We will consider cones spanned by these characteristic
vectors. We therefore agree to denote by $V(\Ss)$ the cone spanned by
the vectors $v_X$ for $X\in \Ss$, $\Ss\,{\subseteq}\,\LL$.

Let~$\LL$ be a finite atomic lattice and~$\GG$ a building set
in~$\LL$. We define a rational, polyhedral fan $\Sigma(\LL,\GG)$ in $\R^n$ 
by taking cones $V(\Ss)$  for any nested set 
$\Ss$ in~$\LL$,
\begin{equation} \label{def_Sigma}
 \Sigma(\LL,\GG) \, \, :=\, \, \{\, V(\Ss)\,|\, \Ss\in \NN(\LL,\GG)\,\} \, . 
\end{equation}  
By definition, rays in $\Sigma(\LL,\GG)$ are in $1$-$1$ correspondence 
with elements in~$\GG$; the face poset of $\Sigma(\LL,\GG)$ coincides with
the face poset of $\NN(\LL,\GG)$. To specify the set of cones
in~$\Sigma(\LL,\GG)$ of a fixed dimension~$k$, or nested sets in~$\GG$
with~$k$ elements, we often use the notation $\Sigma(\LL,\GG)_k$ or
$\NN(\LL,\GG)_k$, respectively.

\begin{proposition} \label{prop_unimodularity}
The polyhedral fan $\Sigma(\LL,\GG)$ is unimodular.
\end{proposition}

\begin{pf} We need to show that for any nested set $\Ss\,{\in}\,\NN(\LL,\GG)$
the set of generating vectors for $V(\Ss)$,
$\{\,v_X\,|\,X\,{\in}\,\Ss\}$, can be extended to a lattice basis
for~$\Z^n$. To that end, fix a linear order $\prec$ on $\Ss$ that
refines the given order on~$\LL$, and write the generating vectors
$v_X$ as rows of a matrix $A$ following this linear order. Now
transform $A$ to a matrix $\widetilde A$, replacing each vector $v_X$
by the characteristic vector $v_{\widetilde X}$ of $\widetilde X$,
with
\[
     \widetilde X \, \,= \, \,  \bigvee_{{Y\in \Ss}\atop{Y\preceq X}}\, Y \,.
\]
For each~$X$ this can be done by adding rows $v_Z$ to $v_X$ for elements
\[
Z\, \in\,  {\rm max}_{\LL}\{\, Y\in \Ss\, | \, 
         Y\,{\prec}\,X,\, Y \mbox{ incomparable to }\, X \mbox{ in }\,\LL\}\, ,
\]
the reason being that characteristic sets of atoms for incomparable elements 
of a nested set are disjoint~\cite[Prop.~2.5(1),\,2.8]{FK}. 
The matrix $\widetilde A$ 
clearly has rows with strictly increasing  support, hence can be easily 
extended to a square matrix with determinant~$\pm1$. The same extra rows will 
complete the rows of the original matrix~$A$ to a lattice basis for~$\Z^n$.
\end{pf}

\begin{remark} \label{rem_unimodularity}
  In section~\ref{sect_blowups} we will give a more 
  constructive description 
  of $\Sigma(\LL,\GG)$, picturing the fan as the result of successive
  stellar subdivisions of faces of the \mbox{$n$-dim}\-en\-sional cone
  spanned by the standard lattice basis for~$\Z^n$ and subsequent
  removal of faces (compare Thm.~\ref{thm_ST}).  From this
  description, unimodality of the fan will follow immediately.
\end{remark}

Let $X_{\Sigma(\LL,\GG)}$ denote the toric variety associated with
$\Sigma(\LL,\GG)$. If there is no risk of confusion, we will
abbreviate notation by using $X_{\Sigma}$ instead.  $X_{\Sigma}$ is a
smooth, non-complete, complex algebraic variety. Crucial for us will
be its stratification by torus orbits $\OO_{\Ss}$, in one-to-one
correspondence with cones $V(\Ss)$ in $\Sigma(\LL,\GG)$, thus with
nested sets $\Ss$ in $\GG$.

The orbit closures $[\OO_{\Ss}]$, $\Ss{\in}\NN(\LL,\GG)_{n-k}$,
generate the Chow groups $A_k(X_{\Sigma})$, $k\,{=}\,0,\ldots,n$. We
describe generators for the groups of relations among the
$[\OO_{\Ss}]$, $\Ss\,{\in}\,\NN(\LL,\GG)_{n-k}$, in $A_k(X_{\Sigma})$
for later reference.  This description is due to Fulton and
Sturmfels~\cite{FS}. We present here a slight adaptation to our present
context.

\begin{proposition} \label{prop_relAk}
{\rm (\cite[2.1]{FS})} The group of relations among generators
$[\OO_{\Ss}]$, $\Ss\,{\in}\,\NN(\LL,\GG)_{n-k}$, for the $k$-th Chow
group $A_k(X_{\Sigma})$, $k{=}0,\ldots,n$, is generated by relations
of the form
\begin{equation}\label{eq_relAk}
r(\TT,b) \,\, \,\, = 
\sum_{{\Ss \supset \TT}\atop {\Ss \in \NN(\LL,\GG)_{n-k}}}
      <b,z_{\Ss,\TT}> \,\, [\OO_{\Ss}] \, ,  
\end{equation}
where $\TT$ runs over all nested sets with $n\,{-}\,k\,{-}\,1$
elements and $b$ over a generating set for the sublattice determined by
$V(\TT)^{\perp}$ in the dual lattice Hom$(\Z^n,\Z)$.  Here,
$z_{\Ss,\TT}$ is a lattice point in $V(\Ss)$ generating the
($1$-dimensional) lattice
{\rm span}$(V(\Ss)\,{\cap}\,\Z^n)/${\rm span}$(V(\TT)\,{\cap}\,\Z^n)$.
\end{proposition}

Since $X_{\Sigma(\LL,\GG)}$ is non-singular, the intersection
product~$\cdot$ makes ${\rm Ch}^*(X_{\Sigma})\,{=}$ 
 $\oplus_{k=0}^n{\rm Ch}^k(X_{\Sigma})$ with 
${\rm Ch}^k(X_{\Sigma})\,{=}\,A_{n-k}(X_{\Sigma})$ into a commutative
graded ring, the {\em Chow ring\/} of $X_{\Sigma(\LL,\GG)}$.

\begin{theorem} \label{thm_ChD}
Let $X_{\Sigma(\LL,\GG)}$ be the toric variety associated with a finite
atomic lattice~$\LL$ and a combinatorial building set $\GG$ in~$\LL$ as described 
above. Then the assignment
$x_G\mapsto [\OO_{\{G\}}]$ for $G\,{\in}\,\GG$, extends to an isomorphism
\[
  D(\LL,\GG)\, \, \cong \, \, {\rm Ch}^*(X_{\Sigma(\LL,\GG)}) \, .
\]    
\end{theorem}

\begin{pf}
Orbit closures $[\OO_{\{G\}}]$ in $X_{\Sigma}$ that  correspond to the
rays $V(\{G\})$ in $\Sigma(\LL,\GG)$ for $G\,{\in}\,\GG$, generate
${\rm Ch}^*(X_{\Sigma})$ multiplicatively, since
\[
   [\OO_{\Ss}]\, \, = \, \, [\OO_{\{G_1\}}] \,\cdot\, \ldots  
                                             \,\cdot\, [\OO_{\{G_k\}}]  
\]
for $\Ss\,{=}\,\{G_1,\ldots, G_k\}\,{\in}\,\NN(\LL,\GG)$, $\cdot$ 
denoting the intersection product (see \cite[p.100]{F}).  

Moreover, relations as in $D(\LL,\GG)$ hold. Indeed,
the intersection products of orbit closures corresponding to rays that do
{\em not\/} span a cone in $\Sigma(\LL,\GG)$ are~$0$~\cite[p.100]{F},
which is exactly the monomial relations (\ref{rel_monD}) for
non-nested index sets in $D(\LL,\GG)$. Relations (\ref{eq_relAk}) in
${\rm Ch}^1(X_{\Sigma})\,{=}\,A_{n-1}(X_\Sigma)$ as described above
coincide with the linear relations~(\ref{rel_linD}) in $D(\LL,\GG)$
\begin{equation}
   r(\emptyset, v_A) \, \,  = \, \, 
            \sum_{G\in \GG} \,<v_A,v_G>\, [\OO_{\{G\}}] 
                      \,\,   = \,\, 
            \sum_{G\geq A} \,\,[\OO_{\{G\}}]\, , \label{eq_relAn}  
\end{equation}
the $v_{A}$, for $A\,{\in}\,\at(\LL)$, forming a basis for the lattice
orthogonal to $V(\emptyset)\,{=}\,0$ in $\Z^n$.

Thus, sending $x_G$ to $[\OO_{\{G\}}]$ for $G\,{\in}\,\GG$, we have 
a surjective ring homomorphism from $D(\LL,\GG)$ to the Chow 
ring of $X_{\Sigma}$. It remains to show that the relations~(\ref{eq_relAk})
in ${\rm Ch}^*(X_{\Sigma})$ follow from relations~(\ref{eq_relAn}) in 
${\rm Ch}^1(X_{\Sigma})$, and from monomials over non-nested index sets 
being zero. 

Let us fix some notation. For $\TT\,{\in}\,\NN(\LL,\GG)$ and  
$X\,{\in}\,\TT$ define
\[
             \Delta_{\TT}(X)\, \, := \, \,
                            \lfloor X \rfloor \, \setminus \,
        \bigcup_{{Y<X}\atop{Y\in \TT}}\,\lfloor Y \rfloor \, , 
\]
the set of atoms that are below~$X$, but not below the join of all $Y$
in $\TT$ that are smaller than~$X$. Observe that
$\Delta_{\TT}(X)\,{\neq}\, \emptyset$ for any $X\,{\in}\,\TT$, since
$\TT$ is nested, and
\[
            \lfloor \bigvee \TT \rfloor 
                         \, \, = \, \, 
            \bigcup_{X\in \TT} \Delta_{\TT}(X)\, .
\]

For $\TT\,{\in}\,\NN(\LL,\GG)_{k-1}$, $k\,{\geq}\,2$, the sublattice 
determined by  
$V(\TT)^{\perp}$ in the dual lattice
is generated by vectors in $\CC_1\,{\cup}\,\CC_2$, where
\begin{eqnarray*}
  \CC_1 & = & \{\, v_{A_i}\,{-}\,v_{A_j}\,| \, A_i, A_j \in   \Delta_{\TT}(X) 
                                    \mbox{ for some }\, X\in \TT\}\, , \\
  \CC_2 & = & \{\, v_{A} \,| \, 
         A \in \at(\LL)\setminus \lfloor \bigvee \TT \rfloor\,\}\, .
\end{eqnarray*}
Observe that $\CC_1\,{\cup}\,\CC_2$ contains $\sum_{X\in
  \TT}(|\Delta_{\TT}(X)|-1)+  |\at(\LL)\,{\setminus}\, \lfloor
\bigvee \TT \rfloor| = |\at(\LL)| - |\TT| = {\rm codim}\, V(\TT)$
linear independent vectors, thus a basis for the sublattice determined 
by~$V(\TT)^\perp$.

For $\TT\,{\in}\,\NN(\LL,\GG)_{k-1}$, $k\,{\geq}\,2$, and
$v_{A_i}\,{-}\,v_{A_j}\,{\in}\,\CC_1$, the relation~(\ref{eq_relAk})
reads as 
\begin{eqnarray*}
\lefteqn{ \hspace{-0.7cm}
r(\TT,v_{A_i}\,{-}\,v_{A_j}) } \\ 
& =  &  \sum_{{\Ss \supset \TT}\atop {\Ss \in \NN(\LL,\GG)_{k}}}
      <v_{A_i}\,{-}\,v_{A_j},z_{\Ss,\TT}> \,\, [\OO_{\Ss}]  \\
& = &  \sum_{{Y\in \GG\setminus \TT}\atop {\TT \cup \{Y\}\in \NN(\LL,\GG) }}
      <v_{A_i}\,{-}\,v_{A_j},v_Y> \,\, [\OO_{\TT\cup \{Y\}}]  \\
& = &  [\OO_{\TT}]  \cdot   
       \big( \sum_{{Y\in \GG\setminus \TT, Y\geq A_i}\atop 
                        {\TT \cup \{Y\}\in \NN(\LL,\GG) }}
               [\OO_{\{Y\}}] \, \, - \, \, 
             \sum_{{Y\in \GG\setminus \TT, Y\geq A_j}\atop 
                        {\TT \cup \{Y\}\in \NN(\LL,\GG) }}
               [\OO_{\{Y\}}]\,\,\, \big)\, .     
\end{eqnarray*}
Monomials over non-nested index sets being zero, we may drop the
condition $\TT\,{\cup}\,\{Y\}\,{\in}\,\NN(\LL,\GG)$ in both sums.
Moreover, if $Y\,{\in}\,\TT$, $Y$ either is larger than both $A_i$ and
$A_j$, or not larger than either of them. Thus, both sums in
$r(\TT,v_{A_i}\,{-}\,v_{A_j})$ are relations of type~(\ref{eq_relAn}),
hence $r(\TT,c)$, $c\in \CC_1$, is a consequence of relations of type
(\ref{rel_monD}) and (\ref{rel_linD}) holding in ${\rm
  Ch}^*(X_{\Sigma})$, as claimed.

For $v_A\,{\in}\,\CC_2$, the reasoning is similar, but easier. Indeed
\begin{eqnarray*}
r(\TT,v_A) & = & 
      \sum_{{\Ss \supset \TT}\atop {\Ss \in \NN(\LL,\GG)_{k}}}
      <v_A,z_{\Ss,\TT}> \,\, [\OO_{\Ss}]  \\
& = &  
      \sum_{{Y\in \GG\setminus \TT}\atop {\TT \cup \{Y\}\in \NN(\LL,\GG) }}
      <v_A,v_Y> \,\, [\OO_{\TT\cup \{Y\}}]  \\
& = &  [\OO_{\TT}]\, \, \cdot \, \, \sum_{Y\geq A}\, [\OO_{\{Y\}}]\, ,
\end{eqnarray*} 
since no $Y\,{\in}\,\TT$ can be larger than~$A$, and again, by monomials
over non-nested sets being zero, the condition  $\TT \cup \{Y\}\in \NN(\LL,\GG)$
can be dropped. This completes our proof.
\end{pf}


\section{A geometric description of $X_{\Sigma(\LL,\GG)}$}
\label{sect_blowups}

The goal of this section is to give a geometric description of the variety
$X_{\Sigma(\LL,\GG)}$. For an arbitrary atomic lattice $\LL$,
we describe the toric variety
$X_{\Sigma(\LL,\GG)}$ as the
result of a sequence of blowups of closed torus orbits and subsequent
removal of a number of open orbits.
We start with
a more constructive description of the fan 
$\Sigma(\LL,\GG)$ as the result of a sequence of stellar subdivisions
and subsequent removal of a number of open cones.

\smallskip
We allow the same setting as for the definition of $\Sigma(\LL,\GG)$
in~(\ref{def_Sigma}). Let $\LL$
be a finite atomic lattice with set of atoms $\at
(\LL)\,{=}\,\{A_1,\ldots,A_n\}$ and $\GG$ a building set in~$\LL$.

\smallskip
\noindent
{\bf Construction of $\Theta(\LL,\GG)$.}
 \newline 
(0) Start with the fan $\Theta_0$ given by the $n$-dimensional cone 
spanned by the coordinate vectors in~$\R^n$ together with all its faces.
\newline 
(1) Choose a linear order
$\succ$ on $\GG$ that is
non-increasing with respect to the original partial order on~$\LL$,
i.e., $G\,{\leq}\,G'$ implies $G'\,{\succeq}\,G$.
Write $\GG\,{=}\,\{G_1\succ G_2\succ\cdots\succ G_t\}.$ 
Construct a fan $\widetilde \Theta(\LL,\GG)$ by successive barycentric 
stellar subdivisions in faces $V(\lf G_i\rf)$ of $\Theta_0$ for 
$i=1,\ldots,t$, introducing in each step a new ray generated by 
the characteristic vector $v_{G_i}$, $i=1,\ldots,t$. \newline
(2) Remove from $\widetilde \Theta(\LL,\GG)$ all (open) cones $V(\TT)$
with index sets of generating vectors~$\TT$ that are not nested in~$\GG$ and 
denote the resulting fan by $\Theta(\LL,\GG)$.

\begin{theorem} \label{thm_ST}
The fan $\Theta(\LL,\GG)$ constructed above coincides with 
the fan $\Sigma(\LL,\GG)$ defined in section~\ref{sect_torv}. 
\end{theorem}

\begin{pf}
  By construction the fans share the same generating vectors. In fact,
  due to the removal of cones in step (2) of the construction above,
  it is enough to show that for any nested set $\Ss\in \NN(\LL,\GG)$
  there exists a cone in $\widetilde \Theta(\LL,\GG)$ containing 
  $V(\Ss)$ as a face. Due to the recursive construction of 
  $\widetilde \Theta(\LL,\GG)$ this statement reduces to the 
  following claim.

\smallskip
\noindent
{\bf Claim.}
Let $\Ss=\{X_1,\ldots,X_k\}$ be nested in $\LL$ with respect to~$\GG$,
and assume that the indexing is compatible with the linear order $\succ$
on~$\GG$, i.e., $X_1\,{\succ}\, \ldots \,{\succ}\,X_k$. For notational convenience,
extend the set by $X_{k+1}\,{:=}\,\hat 0$. 
Then any stellar 
subdivision in $V(\lf G \rf)$, $G\,{\in}\,\GG$, during the construction of 
$\widetilde \Theta(\LL,\GG)$, for $G\,{\succ}\,X_i$, 
$G\,{\not \succeq}\,X_{i-1}$, $i\,{=}\,1,\ldots,k{+}1$, 
retains a cone $W_G$ with
\[
     V(\,\{X_1,\ldots,X_{i{-}1}\}\, \cup \, \lf X_i\rf \, \cup \, \ldots 
                             \, \cup \, \lf X_k\rf\, )
\]
among its faces and for $G\,{=}\,X_i$, $i\,{=}\,1,\ldots,k$,
creates a cone $W_{X_i}$ with
\[ 
     V(\,\{X_1,\ldots,X_i\}\, \cup \, \lf X_{i+1}\rf \, \cup \, \ldots 
                             \, \cup \, \lf X_k\rf\, )
\]
among its faces.

\smallskip
\noindent
{\bf Proof of the claim.} Assume first that $G\,{\succ}\,X_i$, 
$G\,{\not \succeq}\,X_{i-1}$, for some $i{\in}\{ 1,\ldots,k{+}1\}$ 
(the second condition being empty for $k\,{=}\,1$), and assume that 
the previous subdivision step in $V(\lf G' \rf)$, $G'\,{\in}\,\GG$, 
has created, resp.\ retained a cone~$W_{G'}$ with 
$ V(\,\{X_1,\ldots,X_{i{-}1}\}\,{\cup}\, \lf X_i\rf \,{\cup}\,\ldots 
                             $ $ \,{\cup}\, \lf X_k\rf\, )$
among its faces.

If $W(G')$ does not contain $V(\lf G \rf)$, it will
not be altered by stellar subdivision in~$V(\lf G \rf)$. Any cone that
is to be altered when subdividing~$V(\lf G \rf)$ needs to be contained
in star$\,V(\lf G \rf)$, hence among its faces needs to contain~$V(\lf G
\rf)$.

If $W(G')$ does contain $V(\lf G \rf)$ among its faces, choose 
\begin{equation} \label{eq_choiceg}
 g\, \, \in \, \, \lf G \rf\, \setminus \, \bigcup_{j=i}^k\, \lf X_j\rf\, .
\end{equation}
If the set was empty, we would have 
$\lf G \rf\,{\subseteq}\, \bigcup_{j\geq i}\, \lf X_j\rf$, in particular,
\[
   G \, \, \leq \, \, \bigvee_{j\geq i}\, X_j 
     \, \, \leq \, \, \bigvee_{{\rm max}\Ss_{\succeq X_i}}\, X_j\, . 
\]
The join on the right hand side is taken over all $X_j$ that are maximal 
among $X_1,X_2,\ldots,X_i$ with respect to the partial order in~$\LL$. 
Since these elements 
are pairwise incomparable
and nested in~$\LL$ they are the factors
of their join. This implies that
$G\,{\leq}\,X_{j^*}$ for some $j^*\geq i$~\cite[Prop.\,2.5(i)]{FK} contradicting
the fact that $G\succ X_{j^*}$.

Hence we can choose $g$ as described in~(\ref{eq_choiceg}) and, when 
subdividing
$V(\lf G \rf)$, we replace $W_{G'}$ by $W_G$ by substituting the new 
ray~$\la v_G \ra$ for the ray $\la v_g \ra$ in $W_{G'}$. Observe that 
$V(\,\{X_1,\ldots,X_{i{-}1}\}\,{\cup}\,\lf X_i\rf \,{\cup}\,\ldots 
                             $ ${\cup} \, \lf X_k\rf\,)$ 
remains as a face in the newly created cone~$W_G$.

Assume now that $G\,{=}\,X_i$ 
and again denote the cone emerging from the previous subdivision
step by $W_{G'}$, assuming that it contains 
$V(\,\{X_1,\ldots,X_{i{-}1}\}$ $\,{\cup}\,\lf X_i\rf \,{\cup}\,\ldots 
\,{\cup} \, \lf X_k\rf\,)$ among its faces. When subdividing $V(\lf X_i \rf)$
now replace $W_{G'}$ by $W_{X_i}$ by substituting the new ray~$\la v_{X_i} \ra$
for the generating ray associated with some 
\[
    x_i \,\, \in \, \,  
        \lf X_i\rf\, \setminus \bigcup_{j\geq i{+}1}\,\lf X_j\rf 
        \, \, = \, \, 
        \lf X_i\rf\, \setminus \bigcup_{{j\geq i{+}1}\atop {X_j<X_i}}
                                                         \,\lf X_j\rf
        \, \, = \, \,
        \Delta_{\Ss} (X_i)\, , 
\]
where the right hand side is non-empty as we observed before 
(see proof of Thm.~\ref{thm_ChD}).

Note that $V(\,\{X_1,\ldots,X_{i}\}\,{\cup}\,\lf X_{i{+}1}\rf \,{\cup}\,\ldots 
                             \,{\cup} \, \lf X_k\rf\,)$
is a face of the newly created cone~$W_{X_i}$. This completes the proof 
of our claim.   
\end{pf}

\begin{corollary} \label{crl_bu}
The toric variety $X_{\Sigma(\LL,\GG)}$ can be constructed 
as follows.
Start from the toric variety associated with the $n$-dimensional cone 
spanned by the standard lattice basis in~$\Z^n$, i.e., from $\C^n$ 
stratified by torus orbits. Perform a sequence of blowups in orbit 
closures associated with faces $V(\lf G\rf)$ of the standard cone for 
$G\,{\in}\,\GG$ in some linear, non-increasing order. Remove from the 
resulting variety all open torus orbits that correspond to cones in 
$\widetilde\Theta(\LL,\GG)$ indexed with non-nested subsets of~$\LL$.
\end{corollary}

It follows immediately from this description that the toric variety 
$X_{\Sigma(\LL,\GG)}$ is smooth.


\section{Examples} \label{sect_examples}

We discuss a number of examples to illustrate the central notions of this 
article.

\noindent
{\bf Partition lattices.} \newline
Let $\Pi_n$ denote the lattice of set partitions of 
$[n]$ ordered by reversed refinement. As we mentioned above, the
partition lattice $\Pi_n$ occurs as the intersection lattice of the braid 
arrangement $\Aa_{n-1}$ (compare Example~\ref{ex_braid}).

\medskip
\noindent
\begin{minipage}{5.9cm}
For $n\,{=}\,3$, the only building set is the maximal one, i.e., 
$\GG=\Pi_3\setminus\{\hat 0\}$. Denoting elements as in the Hasse 
diagram depicted on the right, the nested set complex~$\NN(\Pi_3,\GG)$ 
contains the following simplices:
\end{minipage} \mbox{ }\hspace{0.6cm} 
\begin{minipage}{3cm}
\vspace{-0.4cm}
  \begin{picture}(0,0)%
   \includegraphics{Pi3.pstex}%
  \end{picture}%
 \input{Pi3.pstex_t}%

\end{minipage}

\[
   \NN(\Pi_3,\GG) = \{H_{12}, H_{13}, H_{23}, U, H_{12}U, 
                           H_{13}U, H_{23}U \}\, .
\]

\medskip
\noindent
The algebra $D(\Pi_3,\GG)$ thus is the following:
\begin{eqnarray*}
\lefteqn{ \hspace{-1cm}
   D(\Pi_3,\GG) \, \, = \, \, 
       \Z\,[x_{H_{12}}, x_{H_{13}}, x_{H_{23}}, x_{U}] \, \Big/ \, } \\
   & &   \qquad \qquad  \qquad \left\langle \begin{array}{l}
        x_{H_{12}}x_{H_{13}},\,\, x_{H_{12}}x_{H_{23}},\,\,  
        x_{H_{13}}x_{H_{23}} 
                                                                            \\
        x_{H_{12}}+x_U,\,\,  x_{H_{13}}+x_U, \,\, x_{H_{23}}+x_U
       \end{array}
       \right\rangle \, . 
\end{eqnarray*}
\noindent
We find that $D(\Pi_3,\GG)\cong \Z\,[x_{U}]\, / \, \langle x_U^2\rangle$, which
illustrates Corollary~\ref{crl_HDP}. The compactification 
$\overline{Y}\!_{\Pi_3\setminus \{\hat 0\}}$ of the complement of the 
projectivized braid arrangement 
$\P\Aa_{2}$ (a three times punctured $\CP^1$) is the complex projective line. 

\smallskip
\noindent
\begin{minipage}{6cm}
To visualize the fan $\Sigma(\Pi_3,\GG)$ we choose to depict its intersection 
with a hyperplane orthogonal to the diagonal ray
in the positive octant of $\R^3$. To shorten notation, we denote rays by 
building set elements. 
\end{minipage}
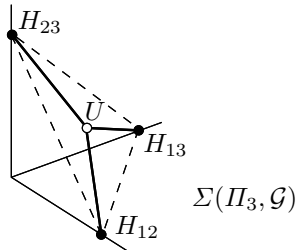
\mbox{ } \hspace{1cm} 
\begin{minipage}{3cm}
  \begin{picture}(0,0)%
   \includegraphics{SPi3.pstex}%
  \end{picture}%
 \input{SPi3.pstex_t}%

\end{minipage}

\noindent
The toric variety $X_{\Sigma(\Pi_3,\GG)}$ is the blowup of $\C^3$ in $0$ 
with open torus orbits corresponding to cones
$V(H_{12}, H_{13})$, $V(H_{12}, H_{23})$, $V(H_{13}, H_{23})$
and $V(H_{12}, H_{13},U)$, $V(H_{12}, H_{23},U)$, $V(H_{13}, H_{23},U)$
removed. What we remove here, in fact, are the proper transforms of the 
three coordinate axes of $\C^3$ after blowup in $0$.

\medskip
\noindent
For $n\,{=}\,4$, we have several choices when fixing a building set. The
partitions with only one non-trivial block of size~$\geq 2$ 
form the minimal building set $\FF$. To obtain the others we add
any number of
$2$-block partitions in $\Pi_4$.

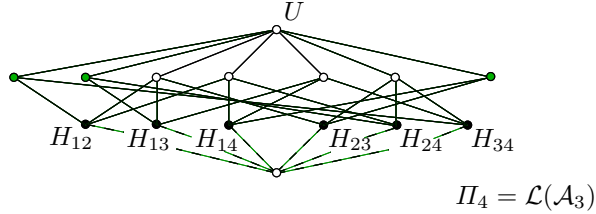
\vspace{0.2cm}
\begin{center}
  \begin{picture}(0,0)%
   \includegraphics{Pi4.pstex}%
  \end{picture}%
 \input{Pi4.pstex_t}%

\end{center}

The nested set complex $\NN(\Pi_4,\FF)$ is a $2$-dimensional complex on $11$ 
vertices. It is a cone with apex $U$, the simplices in its base 
$\NN(\Pi_4,\FF)_0$ being the ordered subsets in $\FF\setminus\{U\}$ together 
with the pairs $H_{12}H_{34}$, $H_{13}H_{24}$, $H_{14}H_{23}$.
We depict below the $1$-dimensional base $\NN(\Pi_4,\FF)_0$.  
To simplify notation we label vertices with the non-trivial block 
of the corresponding partition. The non-ordered nested pairs appear
shaded. 

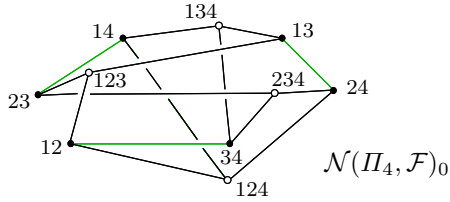
\begin{center}
  \begin{picture}(0,0)%
   \includegraphics{Npi4min.pstex}%
  \end{picture}%
 \input{Npi4min.pstex_t}%

\end{center}

Choosing instead of $\FF$ the maximal building sets $\GG$ in $\Pi_4$, i.e., 
including the $2$-block partitions into the building set, results 
in a subdivision of these edges by additional vertices $H_{12|34}$,  
$H_{13|24}$ and  $H_{14|23}$ corresponding to the newly added building set 
elements.

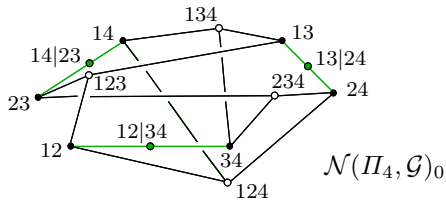
\begin{center}
  \begin{picture}(0,0)%
   \includegraphics{Npi4max.pstex}%
  \end{picture}%
 \input{Npi4max.pstex_t}%

\end{center}

\medskip
Simplifying the presentation of the algebra $D(\Pi_4,\FF)$ given in 
Definition~\ref{def_D} yields
\begin{eqnarray*}
\lefteqn{ 
   D(\Pi_4,\FF) \, \cong \, 
       \Z\,[x_{123}, x_{124}, x_{134}, x_{234}, x_{U}] \, \Big/  }\\
              & &  \qquad \qquad \qquad \qquad
       \left\langle \begin{array}{ll} 
        x_{ijk} \,x_U & \mbox{for all }\, 1{\leq}i{<}j{<}k{\leq}4 \\
        x_{ijk} \, x_{i'j'k'} & \mbox{for all }\, ijk \neq i'j'k' \\
        x_{ijk}^2 + x_U^2 \quad & \mbox{for all }\, 1{\leq}i{<}j{<}k{\leq}4  
       \end{array}
       \right\rangle \, , 
\end{eqnarray*}
where we index generators corresponding to rank~$2$ lattice elements by
the non-trivial block of the respective partitions. The linear basis 
described in Corollary~\ref{crl_basis} 
is given by the 
monomials $x_{123}$, $x_{124}$, $x_{134}$, $x_{234}$, $x_U$, and $x_{U}^2$.

\medskip
For completeness, we state the description of $D(\Pi_n,\FF)$ for general~$n$, 
where $\FF$ 
again denotes the minimal building set, i.e., the set of $1$-block partitions 
in $\Pi_n$.  Having in mind that $D(\Pi_n,\FF)$ 
is isomorphic  to the cohomology of the Deligne-Knudson-Mumford 
compactification $\overline M_{0,n+1}$  of the moduli space of 
$n{+}1$-punctured complex projective lines
(compare Example~\ref{ex_braid}), the following 
presentation should be compared with presentations for 
$H^*(\overline M_{0,n+1})$ given earlier by Keel~\cite{Ke}.

We index generators for $D(\Pi_n,\FF)$ with subsets of $[n]$ of
cardinality larger than two representing the non-trivial blocks in the
respective partitions and obtain:
\begin{eqnarray*}
\lefteqn{ 
  D(\Pi_n,\FF) \,\,  \cong  \,\, 
       \Z\,[\,\{x_S\}_{S\subseteq [n], |S|\geq 2}\,] \,\, \Big/}\\
    & &  \qquad \qquad \qquad \qquad \qquad 
        \left\langle \begin{array}{ll}
        x_S \, x_T & \mbox{for }\, 
                          S\cap T \neq \emptyset, \\
                & \mbox{and }\, 
                          S\not\subseteq T, T\not\subseteq S\, ,    \\[0.1cm]
        \sum_{\{i,j\}\subseteq S} \, x_S \quad & 
                                          \mbox{for }\, 1\leq i<j\leq n\,
       \end{array}
       \right\rangle \, .     
\end{eqnarray*}

\medskip
\noindent
{\bf A non-geometric lattice.}

\noindent
\begin{minipage}{7cm} 
Consider the lattice $\LL$ depicted by its 
Hasse diagram on the right.
We obtain the following building sets:
\begin{eqnarray*}
\GG_1&=& \{A_1, A_2, A_3, U\}\, ,\\
\GG_2&=& \{A_1, A_2, A_3, Y_1,U\}\, ,\\
\GG_3&=& \{A_1, A_2, A_3, Y_1,Y_2,U\}\, ,
\end{eqnarray*}
the only other choice being to replace $Y_1$ by $Y_2$ in~$\GG_2$.
\end{minipage}
\hfill
\begin{minipage}{4cm}
  \begin{picture}(0,0)%
   \includegraphics{non-geom.pstex}%
  \end{picture}%
 \input{non-geom.pstex_t}%

\end{minipage}

\medskip
For a description of the nested set complexes we refer to the 
corresponding fans $\Sigma(\LL,\GG_i)$, $i{=}1,2,3$, shown below. 
The standard presentations for $D(\LL,\GG_i)$, $i{=}1,2,3$, according to
Definition~\ref{def_D}
simplify so as to reveal the Hilbert functions of the algebras to be
\[
H(D(\LL,\GG_i),t)\, = \, 1\,+\, i\, t \qquad \mbox{for }\, i=1,2,3\, ,
\] 
with basis in degree~$1$ being the generators associated to building set 
elements other than atoms.

\smallskip
We depict the fans $\Sigma(\LL,\GG_i)$, $i{=}1,2,3$, 
again by drawing their intersections with a hyperplane orthogonal 
to the diagonal ray in the positive octant of~$\R^3$.

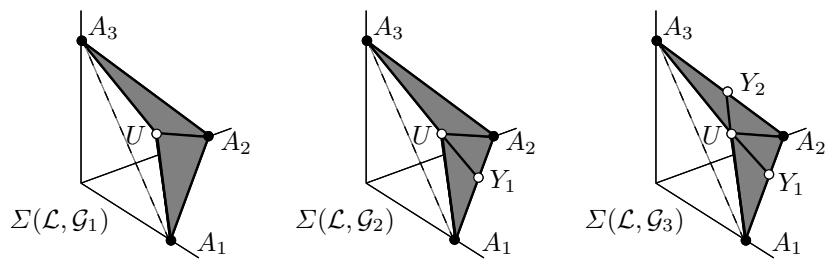
\vspace{0.2cm}
\begin{center}
  \begin{picture}(0,0)%
   \includegraphics{Snon-geom.pstex}%
  \end{picture}%
 \input{Snon-geom.pstex_t}%

\end{center}

The toric variety $X_{\Sigma(\LL,\GG_1)}$ is the result of blowing up
$\C^3$ in the origin, and henceforth removing the open torus orbits 
corresponding to one original $2$-dimensional cone and the unique 
$3$-dimensional cone  containing it. 

The toric varieties $X_{\Sigma(\LL,\GG_2)}$ and $X_{\Sigma(\LL,\GG_3)}$
differ from $X_{\Sigma(\LL,\GG_1)}$ by blowups in one, resp.\ two 
of the original $1$-dimensional torus orbits before removing open orbits 
as above. 


\end{document}